\documentclass[11pt]{amsart}

\title{Symplectic automorphisms of $\mathbf{CP^2}$ and the Thompson group $T$}
\author{Alexandr Usnich}
\address{University Paris VI, 175, rue Chevaleret, Paris, 75013, France}
\email{usnich@ihes.fr}

\begin{document}

\maketitle

 \section{Introduction}
 Study of the group of birational automorphisms of $\mathbf{CP}^2$
preserving the meromorphic $2$-form $\frac{dx\wedge dy}{xy}$ leads
very fast to the understanding of the presence of piecewise linear
geometry behind and also of deep relations with cluster algebras.
Actually, there exists as explained in the paper a morphism from the
group $Symp$ of such birational automorphisms to the group of
piecewise linear automorphisms of $\mathbf{Z}^2$, which is a famous
Thompson group $T$. One particular simple automorphism
$P:(x,y)\mapsto (y,\frac{1+y}{x})$ looks also simply in piecewise
linear world and gives rise to comparisons with mutations. So we
study the subgroup of $Symp$ generated by this $P$ and
$SL(2,\mathbf{Z})$. The main result of this paper is the
construction of a linear representation of this subgroup, the space
where it acts is the inductive limit of Picard groups of the
projective system of surfaces, on which this subgroup acts by
automorphisms. This linear representation is very similar to the
representation of the mapping class group of some cluster algebra,
for which all the seeds are the same and the cluster variables of a
seed are parameterized by $\mathbf{Z}^2$.\\
  For our needs we have to
use the Thompson group $T$, so in the Appendix we find another
representation of it in terms of generators and relations, which
seems to be new and better adopted for a piecewise linear
interpretation.\\
 The paper is organized in the following way: first we construct a
morphism from the group $Symp$ of symplectic automorphisms of
$\mathbf{CP}^2$ to the Thompson group $T$, then we gather some
information about the later, next we find geometrically presentation
of $H$, which is some natural subgroup of $Symp$.\\
 Appendix contains a presentation of $T$ in terms of elements, natural
in our context.\\
 I would like to thank my advisor M.Kontsevich for many useful
 discussions and inspiration.\\

 \subsection{Notations}
 $X$,$Y$,$Z$ homogeneous coordinates on $\mathbf{CP}^2$,
 we will mainly use affine coordinates $x=\frac{X}{Z}$, $y=\frac{Y}{Z}$.\\
 $Symp$ - group of birational automorphisms of $\mathbf{CP}^2$  preserving
 the $2$-form $\omega=\frac{dx\wedge dy}{xy}$, note that $\omega=d(\log{x})\wedge d(\log{y})$.\\
 $P\in Symp$ is defined as $(x,y)\mapsto (y,\frac{1+y}{x})$.\\
 For a form $F$ on $\mathbf{CP}^2$ by $(F)$ we mean it's
 zero locus.\\
 $T$ will denote the Thompson group $T$. It is the group of piecewise
 linear automorphisms of $\mathbf{Z}^2$, for precise definition
 and it's equivalent constructions see one of the chapters below.\\
 $\Gamma=SL(2,\mathbf{Z})$.\\
 $H$ -subgroup of $Symp$ generated by $\Gamma$ and $P$.\\
 $S$ is the subset of $\mathbf{Z}^2$ consisting of vectors with
co-prime coordinates, or in other words
 $S=\mathbf{Z^2}\backslash\bigcup_{n\geq 2}n\mathbf{Z^2}$,
 sometimes we call this set {\it primitive vectors}.\\
 $u\wedge v$ will denote anti-symmetric bilinear product on
 $\mathbf{Z^2}$, normalized in a way, that $(1,0)\wedge(0,1)=1$.
 It also gives an orientation.\\
 For a surface $A$ we denote $Pic(A)$ the Picard group of $A$
 which is a quotient of the group of divisors by the linear
 equivalence.\\
 $\mu$ is a piecewise linear transformation of $\mathbf{Z}^2$,
 defined as $(x,y)\mapsto(x-min(0,y),y)$.\\
 We use standard notation for matrices $C=\begin{pmatrix}-1 & 1 \\ -1 &
 0\end{pmatrix}$,  $I=\begin{pmatrix} 0 & -1 \\ 1 & 0 \end{pmatrix}$, $U=\begin{pmatrix} 1 & 1 \\ 0 & 1 \end{pmatrix}$.\\

\section{Birational automorphisms}

 \subsection{Examples of symplectomorphisms}
 Let us define an action of $P$, $SL(2,\mathbf{Z})$ and
$\lambda\in\mathbf{C}^*$ on $\mathbf{CP}^2$ by birational
automorphisms in the following way:
 $$P:(x,y)\mapsto (y,\frac{1+y}{x})$$
 $$\begin{pmatrix}a&b\\c&d\end{pmatrix}:(x,y)\mapsto (x^a y^b, x^c y^d)$$
 $$\lambda:(x,y)\mapsto(\lambda x,y)$$
 These transformations belong to the group $Symp$.\\
 Using representation of $\omega=d(\log{x})\wedge d(\log{y})$ it
is easy to verify that under these transformations $\omega$ is
indeed preserved, also it is easy to see that any transformation of
the sort $(x,y)\mapsto(Q(y)x,y)$ where $Q$ is any rational function,
may be obtained as a combination of these more primitive ones.\\
 \textbf{Conjecture $1$.}\\
 $Symp$ is generated by transformations of these three types.\\
\\
 To justify this particular form of $P$ we want to observe that
 $P^5=1$\\
 Also we will be interested in the structure of the group $H$
 which is generated by the automorphisms of first two types. If we
 denote by $C$ the matrix $\begin{pmatrix}-1 & 1 \\ 1 &
 0\end{pmatrix}$ and by $I$ the matrix
 with rows $\begin{pmatrix} 0 & -1 \\ 1 & 0 \end{pmatrix}$, than the following relations are easy to verify:
 $$C^3=I^4=[C,I^2]=1$$
 $$PCP=I$$
 $$P^5=1$$
 One recognizes in the first row the defining relations of
 $\Gamma=SL(2,\mathbf{Z})$.\\

 \textbf{Conjecture $2$.}\\
 This list of relations for $H$ is complete, i.e. all the other
relations are consequences of these.\\

 \subsection{Morphism $Symp\rightarrow T$}
 We will provide two ways to see this morphism, which may be
called tropicalization of automorphism. One is the following: take
automorphism
 $$\gamma:x,y\mapsto\frac{P(x,y)}{Q(x,y)},\frac{R(x,y)}{Q(x,y)}$$
 $P,Q,R$ - are polynomials and put $x=t^a$, $y=t^b$, for $t\in(0,+\infty)$, $a,b\in\mathbf{R}$. We are
 interested in the behavior of this automorphism for small generic $t$, so
 we may ignore coefficients and look at exponents, and we may
 search the monomial in $P,R$ with the smallest exponent and
 the monomial in $Q$ with the biggest exponent. These exponents may
 also be found as the smallest powers of $t$ in the development of
 $\gamma$ in series for a small $t$. So we send $(a,b)$ to the
 corresponding exponents and thus obtain a piecewise linear map
 from $\mathbf{R}^2$ to $\mathbf{R}^2$ which is actually defined
 on $\mathbf{Z^2}$. To illustrate this with example, let us take
 $P:x,y\mapsto y, \frac{1+y}{x}$.\\
 $t^a,t^b\mapsto t^b, \frac{1+t^b}{t^a}$, the corresponding
 piecewise linear map is $(a,b)\mapsto (b,min(0,b)-a)$.\\
  Another way to construct this map of groups is more geometric.\\

 \subsection{Geometry of automorphism}
 Here we use standard algebraic geometry techniques to understand
the geometry of a symplectic automorphism. Identifying the curves on
surface with it's strict transform after blow-up, we may speak about
all curves where the form $\omega$ has a pole. All these curves may
be parameterized by the set $S$ of primitive vectors of
$\mathbf{Z}^2$, so the birational symplectomorphism just permutes
these curves, and so defines a piecewise linear automorphism of
$\mathbf{Z}^2$.\\
 Suppose $f\in Symp$, than it is known that $f$ may be presented in a unique way as a
 sequence of blow-ups and then blow-downs, namely there exists $\pi_i:X_i\rightarrow \mathbf{CP}^2$
$i=1,2$ and biregular isomorphism $\phi:X_1\rightarrow X_2$ such
that $\pi_i$ is a composition of blow-ups and we have
$f=\pi_2\circ\phi\circ\pi^{-1}_1$.\\
Denote by $E_1$, $E_2$ exceptional set of $\pi_1$ and $\pi_2$
respectively, so $E_i$ is a union of $\mathbf{P}^1$'s, $\pi_i(E_i)$
is a set of points and $\pi_i$ induces bijection between
$X_i\setminus E_i$ and $\mathbf{CP}^2\setminus \pi_i(E_i)$. Remind
that $f$ is a well defined morphism except at a finite set of
points(precisely at $\pi_1(E_1)$), so we may pull-back the $2$-form
by $f$ everywhere except this finite set and then continue it
analytically to all surface. So $f^*\omega=\omega$ is equivalent to
$\pi^*_1\omega=\phi^*\pi^*_2\omega$ on $X_1\backslash E_1$.\\
The divisor corresponding to the form $\omega$ on $\mathbf{CP}^2$ is $-(X)-(Y)-(Z)$.\\
Let us call {\it chain} a union of smooth rational curves on
surface: $D_1$, $D_2$, \dots ,$D_n$, such that $D_i$ intersects only
with $D_{i+1}$ and $D_{i-1}$ transversally at one point.
$D_{n+1}=D_1$. For example the support of the divisor of $\omega$ on
$\mathbf{CP}^2$ is a chain. We will consider pairs $(X,D)$ where $X$
- surface, $D$ - chain on it. We start from a pair
$(\mathbf{CP}^2,(X)\cup(Y)\cup(Z))$ and when $\widetilde{X}$ is a
blow-up of $X$ at point $p$ and $\widetilde{D}$ is a strict
transform of $D$ then if $p$ is an intersection point of two curves
of a chain, we pass to $(\widetilde{X},\widetilde{D}\cup E)$, $E$
exceptional divisor,
otherwise we pass to $(\widetilde{X},\widetilde{D})$.\\
\textbf{Lemma}\\
 $1$. $\pi_1$ is a sequence of blow-ups which are
supported only on
the chain\\
$2$. chains on $X_1$ and $X_2$ are isomorphic through $f$
\\

\textbf{Proof}\\
First let us make local computation: fix $u$,$v$ local coordinates
and suppose that locally the form looks like $u^k v^l du\wedge dv$.
After blow up this form will become $u^{k+l+1}(\frac{v}{u})^l
du\wedge d(\frac{v}{u})$ which shows that the exceptional divisor
will be a zero of order $k+l+1$.\\
 If we identify curves with their strict transforms, form $w$ will be of order
$-1$ precisely on the curves of chains, so the equality $\pi_1^*
w=\pi_2^* w$ implies the second statement of the lemma.\\
 To prove the first statement, remark that if we blow-up something
else than a point on the chain, then the order of $w$ at the
exceptional divisor becomes positive, so as the positive parts of
divisors of $\pi_1^* w$ and $\pi_2^* w$ are isomorphic we may
contract them simultaneously.\\

\textbf{Corollary.} $\pi_i$ involves blow-ups at two types of
points: first at the points of intersection of curves of the
chain, thus we enlarge chains;\\
second  - at the interior points of curves of the chain.\\

 Denote by $S$ the set of curves of the chain identified with their
strict transforms. Now we realize $S$ as subset of $\mathbf{Z}^2$.
Send curves $(X=0)$, $(Y=0)$, $Z=0$ on $\mathbf{CP}^2$ to points
$(1,0)$, $(0,1)$, $(-1,-1)$ respectively. and use the following
rule: if the curves $C_i$, $C_{i+1}$ on some chain go to
$u,v\in\mathbf{Z}^2$, then send the exceptional curve appearing
after the blow-up of their intersection point to $(u+v)$. This
procedure gives a bijection between curves on chains and the
primitive vectors in $\mathbf{Z^2}$.\\

 Birational symplectic automorphism takes one chain to another, so it takes
corresponding elements in $S$, say $s_1,\dots,s_n$, to some other
elements, preserving cyclic order. But if we blow-up the
intersection point of curves, corresponding to $s_i,s_{i+1}$ then
the exceptional curve(which is corresponding to $s_i+s_{i+1}$) will
go to the exceptional curve of the blow-up of intersection of images
of $s_i$ and $s_{i+1}$. This argument shows that the action on $S$
is linear between $s_i$ and $s_{i+1}$, so we have constructed a
piecewise linear automorphism of $S$ and thus of $\mathbf{Z}^2$,
which is an element of the Thompson group $T$.\\
 In the appendix there is a presentation of $T$ in terms of images of $C$
and $P$. We can check, that $C,I\in SL(2,\mathbf{Z})$ go to
corresponding elements in $T$.\\
 Element $(x,y)\mapsto(y,\frac{1+y}{x})$ maps to $L$.

\section{Thompson group $T$}
 Here we gather some information about group $T$. The classical
 reference is \cite{CFP}.\\

 This group may be interpreted in different ways, as:\\
 - piecewise linear automorphisms of $\mathbf{Z}^2$\\
 - piecewise linear dyadic homeomorphisms of circle\\
 - pairs of binary trees with cyclic bijection between their
 leaves\\
 - pairs of regular trivalent trees with fixed oriented edge and
 cyclic bijection between their leaves\\
 - piecewise projective automorphisms of circle\\
 Here we give five definitions of $T$, sketching the way to pass
from one definition to another.\\

 Let us say that the set $\{s_1,\dots,s_n\}\subset S$ is cyclically
ordered if an hour hand attached to $(0,0)$ will meet them in this
order: $s_1,s_2,s_3,\dots,s_n,s_1,s_2\dots$. They give a
decomposition of $\mathbf{Z}^2$ into cones $\mathbf{Z_{\geq
0}}s_i+\mathbf{Z_{\geq 0}}s_{i+1}$. We require that $s_{i+1}\wedge
s_i=1$ and consider the bijections of $S$ preserving cyclic order
$f$, such that $a s_i+b s_{i+1}$ is mapped to $a f(s_i)+b
f(s_{i+1})$, for $a,b$ - non-negative. We call such bijection a
piecewise linear automorphism
of $\mathbf{Z}^2$.\\

\textbf{Definition1} The Thompson group $T$ is a group of piecewise
linear automorphisms of $\mathbf{Z}^2$.\\

 Dyadic number is a number of the form $\frac{p}{2^q}$,
$p,q\in\mathbf{Z}$. Circle is an interval $[0,1]$ with $0=1$. Dyadic
automorphism of circle would be a piecewise-linear automorphism,
that sends dyadic numbers to dyadic and whose derivatives on the
intervals of linearity are $2^k$, where $k$ - integer numbers.\\

\textbf{Definition2} The Thompson group $T$ is a group of dyadic
automorphisms of a circle.\\

 Actually if we put $S^1\ni 0=1\mapsto (1,0)\in S$,
$\frac{1}{2}\mapsto (0,1)$, $\frac{3}{4}\mapsto (-1,-1)$, then we
can extend this into bijection between dyadic points and $S$
preserving cyclic ordering. Then it is possible to check that the
automorphisms of the set of dyadic points, coming from dyadic
automorphisms, and automorphisms of $S$, coming from PL
automorphisms of $\mathbf{Z}^2$, coincide.\\

 Binary tree has a root, and every vertex has either no or two
descendants - left and right. Vertices with no descendants are
called leaves, they are naturally cyclically ordered. Such a tree
gives a decomposition of the interval: we attach $[0,1]$ to the
root. If vertex has an interval $[p,q]$ attached, we associate
$[p,\frac{p+q}{2}]$ to it's left descendant and $[\frac{p+q}{2},q]$
to it's right descendant. Clearly enough all such intervals are of
length $\frac{1}{2^k}$ and intervals attached to leaves give a decomposition of $[0,1]$.\\

\textbf{Definition3} The Thompson group $T$ has elements represented
by equivalence classes of pairs of binary trees $(R,U)$ with the
equal number of leaves, and a bijection between the leaves of $R$
and the leaves of $U$, preserving the cyclic order. Attaching
simultaneously two descendants to the leaf on $R$ and to the
corresponding by the bijection leaf of $U$ gives an equivalent pair
of trees, so it would give the same element of $T$. A composition of
elements $(R,S)$ and $(M,N)$ is the following: find a binary tree
$Z$, which contains all the vertices of $S$ and $M$, then find
representatives of the equivalence classes of $(R,S)$ and $(M,N)$,
which look like $(X,Z)$ and $(Z,Y)$, by adding descendants to
appropriate leaves, completing $S$ and $M$ to $Z$. The composition
of $(X,Z)$ and $(Z,Y)$ is $(X,Y)$, with the bijection on leaves induced.\\
 If we attach intervals to a tree $(R,U)$ as explained before, then mapping
intervals of leaves of $R$ to intervals of corresponding leaves of
$U$ gives a dyadic automorphism of a circle.\\
\\
 Looking at the binary tree as at the graph, let us remove a root
and it's two adjacent edges and join root's left and right
descendant by an oriented edge going from left to right. All other
edges stay unoriented.\\
 Let $\Upsilon$ be a set of connected planar trees with vertices of valence
$3$ or $1$ and one oriented edge (vertices of valence one are called
leaves and they are cyclically ordered). Consider the set $\Theta$
of pairs $(A,B)$, $A,B\in\Upsilon$ with a bijective identification
of leaves, preserving the cyclic order. Put an equivalence relation
on $\Theta$ generated by the following: if we add two leaves to the
leaf of $A$ and two leaves to the corresponding leaf of $B$, the new
pair of trees with naturally induced bijection of leaves is
equivalent to $(A,B)$.\\

 \textbf{Definition4} The set of equivalence classes of $\Theta$ under this equivalence relation and
composition given on representatives $(A,B)\cdot (B,C)\mapsto (A,C)$
has a well-defined structure of group, and is called the Thompson
group $T$.\\

 Given an element $A\in \Upsilon$, we will attach to it a union of
triangles in the unit disc $D=\{|z|\leq 1\}$. So to every vertex of
$A$ we attach a triangle, and for any two vertices joined by an edge
the corresponding triangles have a common edge. The edges of all
such triangles will be geodesic and the vertices will be rational in
the upper half-plane model of a disc.\\

 To the vertex from which the oriented edge starts we associate
the triangle with vertices $-i,i,-1$, to the vertex where the
oriented edge goes we associate the triangle with vertices $-i,1,i$.
Next, suppose that we associate the triangle $abc$ with a vertex,
and we want to construct a triangle adjacent to the edge $bc$. We
just need to choose a unique hyperbolic automorphism of the disc,
that sends $-1,-i,i$ to $a,b,c$ in this order, then $c,b$ and the
image of $1$ will give the vertices of the new triangle. Triangles
corresponding to leaves will have internal edges - those that
correspond to the only edge of the leaf and by which this triangle
is attached to previously constructed ones. Internal edges form a
polygon, inscribed in the disc, so if we remove the polygon, the
rest is the union of half discs, each arc having a point. So a pair
of trees, corresponding to an element of $\Theta$, gives a pair of
sequences of half discs, thus a
piecewise projective automorphism of a circle.\\
\textbf{Definition5} Thompson group $T$ is a group of piecewise
projective automorphisms of a circle.\\
 Although we consider just finite binary trees, infinite full
binary tree will correspond to the Farey trianguation of a disc.

\section{Presentation in the Picard group}
 The aim of this section is to construct a linear representation of $H$ in an
infinite dimensional vector space. First we construct a pro-scheme
$X$, which is a projective limit of surfaces, blow-ups of
$\mathbf{CP}^2$ in appropriate points. Group $H$ generated by
$\Gamma$ and $P$ will act on it by automorphisms. We make sence of
the Picard group of $X$, and $H$ will act on this free abelian
group, preserving a bilinear symmetric product of signature $(1,\infty)$.\\
 $\Gamma=SL(2,\mathbf{Z})$ is a subgroup of $H$, so first we apply
the same procedure to it and obtain object $Y$, which may be called
universal toric variety.\\
 At the end we just give simple description of action of $\Gamma$
and $\mu=I^{-1}P^{-1}$ on space $V$, which gives a representation of
$H$.\\

\subsection{Universal toric variety}
 In this section we will develop the necessary details about the
projective limit of all toric surfaces. Actually, denote by $Y$ the
pre-scheme, which is a projective limit of blow-ups of
$\mathbf{CP}^2$ at the points of intersection of chains. The group
$\Gamma$ will act on it as the group of automorphisms. Then we
compute the Picard group and find the cones of ample and effective
divisors.\\

 Actually if $f:A\rightarrow B$ is the blow-up of a smooth surface
$B$ at a point then we have the pullback $f^*:Pic(B)\rightarrow
Pic(A)$ and if $(E)$ is the class of the exceptional divisor we have
$Pic(A)=f^*Pic(B)\oplus\mathbf{Z}(E)$ and we have also the
projection morphism $f_*:Pic(A)\rightarrow Pic(B)$ which is just a
forgetting of the exceptional set. So we are allowed to consider an
injective and a projective limit of Pic groups, note them $Pic$ and
$\widehat{Pic}$ respectively, moreover the former is embedded in
later.\\
 Recall that $S$ is the set of points in $\mathbf{Z^2}$ with co-prime coordinates. Actually
$S$ parameterizes the rational curves coming from chains. Let
$PL(S)$ be the space of piecewise linear function from $S$ to
$\mathbf{Z}$, and let $Fun(S)$ be the space of all functions from
$S$ to $\mathbf{Z}$. Inside $PL$ and $Fun$ we have two-dimensional
space $Lin$ of linear functions.\\
 \textbf{Lemma} $Pic(Y)=PL/Lin$, $\widehat{Pic}=Fun/Lin$\\
To explain, why this $Lin$ appears, note that in $Pic(CP^2)$ three
divisors $(X=0)$, $(Y=0)$, $(Z=0)$ are equal. The corresponding
elements in $Pic(Y)$ are functions, generated by there values in the
points $(1,0)$, $(0,1)$, $(-1,-1)$ and extended by linearity to the
rest of the $S$. So $(X=0)$ has corresponding values $1$,$0$,$0$.
$(Y=0)$ has $0$,$1$,$0$, and $(Z=0)$ has $0$,$0$,$1$. So in $Pic(Y)$
these three functions are equal, which is equivalent
for $Lin$ to be trivial.\\

The action $SL(2,\mathbf{Z})$ on $Pic(Y)$ and on $\widehat{Pic}(Y)$
is induced from natural action on $PL$ and $Fun$, space of linear
functions and $Lin$ is preserved.\\
 There is a non-degenerate pairing
 $Pic(Y)\times\widehat{Pic}(Y)\rightarrow\mathbf{Z}$ which would be
defined later. Also there is an embedding
$Pic(Y)\subset\widehat{Pic}(Y)$\\
 The important part of the structure is the cone of effective
divisors $Eff_Y\subset\widehat{Pic}(Y)$ generated by positive
functions on $S$. The dual to it is the cone of ample divisors
$Amp_Y\subset Pic(Y)$ which consists of convex functions. Both
cones are preserved under the action of the group $\Gamma$.\\
 Let $F$ be a piecewise linear function from $S$
to $\mathbf{Z}$. For $a\in S$ we will introduce an index $d(F,a)$
which will measure how much $F$ is far from being linear. So in the
neighborhood of $a$ choose $u$, $v$ such that $u\wedge a=a\wedge
v=1$(so $a$ is between $u$ and $v$). We ask that $F$ is not linear
at most at $a$ on this interval, which can be realized by adding $a$
to $u$ and $v$ sufficiently many times. Then obviously $u+v=ka$ for
$k\in\mathbf{Z}$. Put $d(F,a):=F(u)+F(v)-k F(a)$ and call this value
an index of the function $F$ at $a$. Easy to check that it does not
depend on the choice of $u$, $v$, and it equals
to $0$ at all points where $F$ is linear. $d$ is linear in the first argument.\\

Let us describe a pairing between $Pic(Y)$ and $\widehat{Pic}(Y)$.
Suppose $F$ is a piecewise linear and $G$ is any function on $S$.
The pairing is defined as follows:
$$<F.G>=\sum_{\alpha\in S}d(F,\alpha)G(\alpha)$$
 The reason to introduce $Y$ is that now $\Gamma$ acts by
automorphisms on it. Next we find a pre-scheme where $H$ acts by
automorphisms.\\

\subsection{$X$}

 Notice that when we take a finite sequence of blow-ups of $\mathbf{CP^2}$
at chain points, we come out with a canonical coordinate on each
exceptional divisor. First $X/Y$, $Y/Z$, $Z/X$ are coordinates on
divisors $Z=0$, $X=0$, $Y=0$. Then if we blow-up a point of
intersection of two rational curves with canonical coordinates say
$a$, $b$ on curves $\alpha$, $\beta$ at a point $a=0,b=\infty$,
choose local coordinates $A,B$ at this point, such that
$A\mid_{\alpha}=a$, $A\mid_{\beta}=0$, $B\mid_{\alpha}=0$,
$B\mid_{\beta}=1/b$. Then we take $A/B$ as
canonical coordinate on the exceptional divisor. \\
\textbf{Proposition} $H$ preserves canonical coordinates
on the exceptional divisors coming from chains.\\
\textbf{Proof} It is sufficient to prove this statement for $P,C$
the generators of $H$, which may be done by direct computation.\\
 The corollary is that every rational chain curve in $Y$ has
canonical coordinate on it, so we may consider the set of points
$\Omega$ with coordinate $-1$, and this set is preserved by $H$. The
automorphism $(x,y)\mapsto (y,\frac{1+y}{x})$ involved the blow-ups
of points in this set. So let us consider projective limit $X$ of
surfaces, obtained by blow-ups at $\Omega$. After a blow-up at a
point $p\in\Omega$ replace $p$ by intersection point of divisor of
chain on which $p$ was supported and the exceptional divisor. Then
one may blow-up again in new $p$ and so on. If $s\in S$ corresponds
to the divisor supporting $p$, denote $\delta_s$ the exceptional
divisor. Denote $\delta_s^2$ the exceptional divisor of the blow-up
at the intersection point of $s$ and $\delta_s$ etc. Denote by $X$
the projective limit of blow-ups of this kind and actually this is
the object where $H$ acts by automorphisms.\\

 Let $\sum_{s,k}\mathbf{Z}\delta_s^k$ be a free abelian group
generated by $\delta_s^k$ and $\prod_{s,k}\mathbf{Z}\delta_s^k$
be a free product over indices $s\in S$ and $k\in\mathbf{Z_{>0}}$.\\

 \textbf{Proposition} $Pic(X)=Pic(Y)\oplus\sum_{s,k}\mathbf{Z}\delta_s^k$.\\
$\widehat{Pic}(X)=\widehat{Pic}(Y)\oplus\prod_{s,k}\mathbf{Z}\delta_s^k$.\\

 The pairing between $Pic(X)$ and $\widehat{Pic}(X)$ is extended from that on $Y$ by
saying that $\delta_s^k$ are orthogonal to each other and to
$Pic(X)$ and $\delta_s^k.\delta_s^k=-1$.\\

 Now group $H$ generated by $L$ and $C$ is acting by orthogonal
transformations on $Pic(X)$ ($L$ denotes the inverse of $P$). We
just find it more convenient to describe the action of $L$ on
$Pic(X)$.\\

 Remind that $L$ is given by $(x,y)\mapsto (\frac{1+x}{y},x)$
which is the inverse of $P:(x,y)\mapsto (y,\frac{1+y}{x})$
considered before.\\

 Denote by $A$ the piecewise linear function $A(x,y)=max(0,-y)$.
  Actually it is not linear only at $(1,0)$ and $(-1,0)$.\\
\textbf{Lemma}\\
$L$ acts as follows:\\
$\delta_s^k\mapsto \delta_{L(s)}^k$, for $k\geq 1$, $s\neq (0,\pm 1)$\\
$\delta_{(0,-1)}^k\mapsto\delta_{(1,0)}^{k+1}$\\
$\delta_{(0,1)}^k\mapsto\delta_{(-1,0)}^{k-1}$ for $k\geq 2$\\
$\delta_{(0,1)}^1\mapsto -\delta_{(1,0)}^1+A$\\
for $F\in PL(S)$ we have $F\mapsto F\circ
L^{-1}-F((0,-1))\delta_{(1,0)}^1+F((0,1))(-\delta_{(1,0)}^1+A)$\\

\subsection{The presentation of $H$}

 In $\widehat{Pic}(Y)\subset \widehat{Pic}(X)$ we may consider functions which take value
$1$ at some point and $0$ at all the others. Geometrically they will
define some curve of the chain. As the group $H$ preserve curves of
the chain, the subspace of associated functions will also be
preserved by it. So $H$ will act on the space $V$ dual to it. To
describe it, let us do the following. For any piecewise linear
function $F\in PL$ we associate  $F':=F-\sum_i
d(F,s_i)\delta_{s_i}^1$.\\
\textbf{Lemma} The space $V$ dual to the space of chain curves is
generated by $F'$ and $\delta_{\alpha}^k-\delta_{\alpha}^{k+1}$\\

So $V\subset Pic(X)$ is invariant under group $H$, which may be also
checked directly. the action looks quite simple in the appropriate
basis. First let $F'=F-\sum_{\alpha}d(\alpha,F)\delta_{\alpha}^1$
for any piecewise linear function $F\in PL$. Denote also
$e_{\alpha}^k=\delta_{\alpha}^k-\delta_{\alpha}^{k+1}$. Then such
functions $F'$ modulo linear functions are encoded by their indexes
of non-linearity. So introduce the free abelian groups
$B=\bigoplus_{\alpha\in S}b_{\alpha}$ and
$E=\bigoplus_{k=1,\alpha\in S}^{\infty}e_{\alpha}^k$. Now to any
$F'$ we associate $\sum d(\alpha,F)b_{\alpha}\in B$. Such elements
generate subgroup $B_0\subset B$ of corank $2$. Actually $B_0$
consists of $\sum n_{\alpha}b_{\alpha}$ which are in the kernel of
the map $b\rightarrow\mathbf{Z^2}$ defined by
$b_{\alpha}\mapsto\alpha$ . This is the consequence of the fact,
that pairings of function $F$ with linear functions are $0$.\\
So $H$ acts on the $V=B_0\oplus E$.\\
The action of $\Gamma$ is natural on the indices of $e$ and $b$. The
action of the element $\mu=I^{-1}L$ may be called the mutation
applied to the vector $v=(0,1)$. By analogy we would say, that
$\gamma\circ\mu\circ\gamma^{-1}$ is the mutation applied to the
vector $v$ if $v=\gamma((0,1))$. Here is how the mutation $\mu$
applied to the vector $v\in S$ acts:
$$b_v\mapsto -b_{-v}$$
$$b_{-v}\mapsto e_{-v}+b_{-v}$$
$$b_w\mapsto b_{\mu_v(w)}\mbox{ if }w\wedge v>0$$
$$b_w\mapsto b_{\mu_v(w)}+(v\wedge w)b_{-v}\mbox{ if }w\wedge v<0$$
$$e_v\mapsto b_v+b_{-v}$$
$$e_{v}^k\mapsto e_{v}^{k-1}$$
$$e_{-v}^k\mapsto e_{-v}^{k+1}$$
$$e_w\mapsto e_{\mu_v(w)}$$\\

 Let us make a change of basis:
 if $v\in S$ and $k$ are positive integer then let
 $$p_{kv}:=e_{(k-1)v}+2e_{(k-2)v}+...+(k-1)e_v+k b_v$$
 Now $V$ can be described as the kernel of the map $\bigoplus_{v\in
 \mathbf{Z^2\setminus\{(0,0)\}}}\mathbf{Z}p_v\rightarrow\mathbf{Z^2}$,
 where $p_v$ maps to $v$.\\
 The action of $\mu_v$ is the following: first denote
 $A_v(w):=(w\wedge v)$ if $(w\wedge v)>0$, $0$ if $(w\wedge v)<0$ and $-1$ if it's $=0$\\
 Then $\mu_v(p_w):=p_{\mu_v(w)}+A_v(w)p_{-v}$.\\
 Note that $\mu_v(kv):=(k-1)v$ and $p_0=0$, so $\mu_v(p_{kv}):=p_{(k-1)v}
 -p_{-v}$.\\
 Here we summarize the previous calculations in more compact
formulas to obtain a presentation of $H$ in $V$.\\

 Let $\Lambda=\mathbf{Z^2}\setminus(0,0)$ be a set , and $W=\bigoplus_{i\in
I}\mathbf{Z}e_i$ be a $\mathbf{Z}$-module. Remind the action of
$\mu$ on $\Lambda$ is given by the following: $\mu: i=(x,y)\mapsto (x-min(0,y),y)$.\\
Now we define an action of $\mu$ on $W[q]$ - free abelian group of
polynomials in $q$ with coefficients in $W$.
$$e_{(x,y)}\mapsto e_{(x,y)}+(1-q)y e_{(-1,0)}\mbox{      if }y>0$$
$$e_{(x,y)}\mapsto e_{(x-y,y)}+q y e_{(-1,0)}\mbox{      if }y<0$$
$$e_{(x,0)}\mapsto e_{(x-1,0)}-e_{(-1,0)}$$
 We have a map $W\rightarrow\mathbf{Z^2}$,  $e_{(x,y)}\mapsto(x,y)$,
denote by $V$ it's kernel.\\
\textbf{Lemma}\\
$V$ is an invariant space for the action of $\mu$ and $\Gamma$ and is a presentation of $H$\\
Consider the pull-back of the wedge product from $\mathbf{Z^2}$ to
$W$, and denote it by $\wedge$ too.\\
\textbf{Lemma}\\
$\wedge$ is preserved by the action of $\mu$ and $\Gamma$.\\

\subsection{Cluster mutations}
 By calling a piecewise linear transformation $\mu$ mutation, we made
an allusion to a cluster mutation, so let us try to justify this
similarity.\\
 To a given cluster $\underline{t}$ we associate a vector space
$V(\underline{t})=\bigoplus_{\lambda\in\Lambda(\underline{t})}\mathbf{C}e_{\lambda}$,
where $\Lambda(\underline{t})$ is the set of indices for cluster
variables at this cluster, we suppose $b_{ij}(\underline{t})$ is an
antisymmetric matrix of exponents. Suppose a mutation $\mu_i$ at
variable $i\in\Lambda(\underline{t})$ relates a cluster
$\underline{t}$ to $\underline{t'}$, then we have a bijection
$\sigma_i:\Lambda(\underline{t})\rightarrow\Lambda(\underline{t'})$.
Let us define a map $\mu_i:V(\underline{t'})\rightarrow
V(\underline{t})$ in the following way:\\
$$e_{\sigma(k)}(\underline{t'})\mapsto e_k+max(b_{ik}(\underline{t}),0)e_i$$
$$e_{\sigma(i)(\underline{t'})}\mapsto -e_i$$
It is easy to check that $b(\underline{t})$ which may be considered
as an element of $\bigwedge^2 V^*(\underline{t})$ taking value
$b_{ij}$ at $e_i\wedge e_j$ transforms in a way, consistent with the
transformation of mutation matrix.\\
 Let us now set $\Lambda=\mathbf{Z}^2$ and let us allow
mutations only at the primitive vectors $S\subset\Lambda$. Take all
the seeds to be isomorphic and set $b_{ij}=i\wedge j$, so rather
than having different seeds we would act by mutations on a seed
data. First of all $\sigma_v$ for $v\in S$ acts on $\Lambda$ as
follows:
$\sigma_v:w\mapsto w-min(v\wedge w,0)v$   if $w\wedge v\neq 0$\\
$\sigma_v:w\mapsto w-v$  if $w=kv$\\
It is an automorphism of $\Lambda$. It's easy to see that
$\sigma_{\gamma(v)}=\gamma\circ\sigma_v\circ\gamma^{-1}$ for any
$\gamma\in\Gamma$. We may also denote $\mu:=\sigma_{(1,0)}$.\\
 Recall the alternative description of the Thompson group $T$ given in the
Appendix: $T=<L,C'|I=LC'L, C'^3=I'^4=L^5=(C'I'L)^7=1,
I'^2C'=C'I'^2>$. We put dashes to distinguish these elements from
generators of $SL(2,\mathbf{Z})$.\\
 Let us invert this presentation, i.e. put
$P=L^{-1}$, $C=C'^{-1}$, $I=I'^{-1}$. Then we have:
$$T=<P,C|I=PCP, C^3=I^4=P^5=(PIC)^7=1, I^2C=CI^2>$$
Now we can send the group $\Gamma=SL(2,\mathbf{Z})$ to $T$ by
mapping $\begin{pmatrix}-1&1\\-1&0\end{pmatrix}$ to $C$, and
$\begin{pmatrix}0&-1\\1&0\end{pmatrix}$ to $I$. Denote
$U=CI=\begin{pmatrix}1&1\\0&1\end{pmatrix}$ and $\mu=IP$. In this
notations we deduce the following presentation of $T$:\\

\textbf{Theorem} $T$ is generated by $\Gamma$ and $\mu$, subject to following relations:\\
$(I^2\mu)^2=U^{-1}$\\
$(I^{-1}\mu)^5=1$\\
$(I\mu)^7=1$\\
all the other relations are the consequences of these.\\
 Now the vector space $V_0$ associated to the cluster as explained
above is $V_0=\bigoplus_{s\in\Lambda}\mathbf{C}e_s$, and the action
of $\mu_v$ is:\\
$$e_w\mapsto e_{\sigma_v(w)}+max(v\wedge w,0)e_{-v}$$
$$e_v\mapsto -e_{-v}$$
 And this is the formula for the action of $\mu$ in the representation
of $H$, except for $e_w$ when $w$ is collinear with $v$.\\

\section{Quantization}
 Suppose now that variables $x,y$ are not commuting. Introduce
formal parameter $q$ commuting with both $x$ and $y$ and let $xy=qyx$.\\
Consider the following transformations:
$$P:x,y\mapsto y, qx^{-1}(1+y)$$
$$C:x,y\mapsto x^{-1}yq, x^{-1}q$$
$$I:x,y\mapsto y^{-1}q, x$$
Then, supposing the conjecture $2$, it is easy to verify that all known
relations in group $H$ hold true.\\
 One way to interpret this is the following abstract idea,
inspired by mirror symmetry: we have the family of non-commutative
$P^2$'s, parameterized by $q$, which for $q=1$ degenerate to usual
$\mathbf{CP}^2$ with Poisson bracket, indicating the direction of
noncommutative $P^2$'s. The general framework $[KS]$ predicts that
there is a piecewise-linear object which may be associated to a
degeneration, such that the automorphisms of the family act as
piecewise linear automorphisms of this limit object. In our case,
apparently the limit object is a circle with the set of points $S$
and a piecewise linear structure coming from $\mathbf{Z^2}$. So the
noncommutative automorphisms $P,C,I$ act on this $S^1$ in a
piecewise-linear way.\\

\section{Appendix.  Alternative description of the Thompson group $T$}
 First recall the standard description of $T$ group, given by
generators:
$$A,B,C$$
notations: $R=A^{-1}CB$, $X_2=A^{-1}BA$, $P=A^{-1}RB$\\
 and relations
$$R=B^{-1}C\mbox{   (1)}$$
$$RX_2=BP\mbox{   (2)}$$
$$CA=R^2\mbox{   (3)}$$
$$C^3=1\mbox{   (4)}$$
$$BA^{-1}\mbox{ commutes with }X_2\mbox{   (5)}$$
$$BA^{-1}\mbox{ commutes with }A^{-1}X_2A\mbox{   (6)}$$
as consequences we have the following relations:\\
$R^4=1$, $P^5=1$, $(P^2X_2^{-1})^3=1$, $(PX_2)^4=1$,
$(X_2^2P^{-2})^4=1$.\\
  In what follows we will use the following
procedure: suppose we are given two groups $G_1$, $G_2$, which are
isomorphic and defined by generators and relations, with $N_1$,
$N_2$ - normal subgroups, such that isomorphism takes $N_1$ to
$N_2$. Then we will add some element to $N_1$ and the conjugate of
it's image by isomorphism to $N_2$. We call elements of $N_i$
relations, so that at the end we obtain isomorphic groups $G_1/N_1$
and $G_2/N_2$ which give two presentations of the same group in
terms of generators and relations. If one of the relations has the
generator from one side, and the expression, doesn't depending on
this generator from the other side, we will call it a notation, and
will allow ourselves to remove this generator and this relation and
replace in all the other relations this generator by this
expression. We add the relations to $N_1$ in order to obtain as
$G_1/N_1$ the known representation of $T$ and $G_2/N_2$ will give
another representation of it.\\
 Now we'll search for the representation of $T$-group in terms of
$R$ and $C$. So start from groups $G_1=<A,B,C|R=B^{-1}C, CA=R^2>$
and $G_2=<R,C>$ which are isomorphic, normal groups $N_i$ take
trivial. Thus $A=C^{-1}R^2$, $B=CR^{-1}$. Add $C^3=1$ to both $N_1$
and $N_2$. Adding $R^{-1}A^{-1}CB$ to $N_1$ is equivalent to adding
$R^4$ to $N_2$. $X_2=A^{-1}BA=R^{-2}CCR^{-1}C^{-1}R^2$,
$P=A^{-1}RB=R^{-2}CRCR^{-1}$, so adding $RX_2(BP)^{-1}$ to $N_1$
using $C^3=R^{4}=1$ becomes equivalent to adding $(CR)^{5}$ to
$N_2$. So the group generated by $A,B,C$ with relations $(1)-(4)$ is
isomorphic to the $<R,C|C^3=R^4=(CR)^5=1>$.\\
 Now introduce $\alpha=BA^{-1}=CR^{-1}R^{-2}C=CRC$.
  Again $X_2=A^{-1}BA=R^2(CRC)^{-1}R^2=(R^2\alpha R^2)^{-1}$, $A^{-1}X_2A=R^2CR^2\alpha^{-1}R^2C^{-1}R^2$, so we obtain
another description of $T$-group:\\
generators
$$R,C$$
Notations:$\alpha=CRC$, $\beta=R^2C^{-1}R^2$\\
Then relations are:
$$R^4=1$$
$$(RC)^5=1$$
$$C^3=1$$
$$\alpha\mbox{ commutes with }(R^2\alpha R^2)$$
$$\alpha\mbox{ commutes with }(\beta^{-1}\alpha\beta)$$
\\
 Consider now the element $L=APA^{-1}\in T$.\\
\textbf{Theorem}
$$T=<L,C|I=LCL, C^3=I^4=L^5=(CIL)^7=1, I^2C=CI^2>$$
\textbf{Proof} \\
 Again as in previous considerations let us start from isomorphic
groups $G_1=<R,C|C^3=(RC)^5=1>$ and $G_2=<L,C|C^3=L^5=1>$. The
isomorphism is established by saying that $L\mapsto (RC)^2$ and
$R\mapsto L^{-2}C^{-1}$\\
 Make a notation $I=LCL$ in $G_2$. $R=FC^{-1}=L^{-2}C^{-1}=(CL^2)^{-1}$ so $R^4=1$
becomes $(CL^2)^{-4}=1$ or by conjugating
$1=L(CL^2)^4L^{-1}=(LCL)^4=I^4=1$. So far we have established the
following group isomorphism:
$$<R,C|C^3=R^4=(RC)^5=1>=<L,C|C^3=L^5=(LCL)^4=1>$$
      $L\mapsto (RC)^2$, $C\mapsto C$, $R\mapsto L^{-2}C^{-1}$ \\
 Now consider first commutator(we use $[X,Y]=X^{-1}Y^{-1}XY$, $X^Y=Y^{-1}XY$,
 also by definition of $I=LCL$ we have $LCI=ICL$):\\
  $R^2\alpha R^2=L^{-2}C^{-1}L^{-2}C^{-1}CL^{-2}L^{-2}C^{-1}L^{-2}C^{-1}=
=L^{-2}C^{-1}L^{-1}C^{-1}L^{-2}C^{-1}=L^{-1}I^{-1}C^{-1}L(CL^{-2})^{-1}$,
so $[\alpha,R^2\alpha R^2]\mapsto
[CL^{-2},L^{-1}I^{-1}C^{-1}L(CL^{-2})^{-1}]=[L^{-1}I^{-1}C^{-1}L,(CL^{-2})^{-1}]=
[I^{-1}C^{-1},(IL)^{-1}]^L=(CIILI^{-1}C^{-1}L^{-1}I^{-1})^L=(CI^2L(LCI)^{-1}I^{-1})^L=
(CI^2L(ICL)^{-1}I^{-1})^L=[C^{-1},I^{-2}]^L$, so the first
commutator relation in $G_1$ is equivalent to $[C,I^2]=1$ in
$G_2$.\\
 Note now, that the elements $C,I$ have the relations $C^3=I^4=1,
[C,I^2]=1$, actually they generate the $SL(2,\mathbf{Z})$ inside
$T$, where $C$ corresponds to the matrix with rows $(-1,1),(-1,0)$,
and $I$ - to $(0,-1),(1,0)$.\\
 Let us deal now with the second commutator, we are interested in it
up to conjugation, so we may now conjugate it's arguments.
$R^2=R^{-2}=CL^2CL^2$, $\beta=R^2C^{-1}R^2=CL^2CL^{-1}CL^2$,
$\beta^{-1}\alpha\beta=L^{-2}C^{-1}LC^{-1}L^{-2}C^{-1}CL^{-2}CL^2CL^{-1}CL^2=L^{-1}(LCL)^{-1}L^2C^{-1}LCL^2CLL^{-3}LCLL=
L^{-1}I^{-1}L^2C^{-1}I^2L^{2}IL=(L^2C^{-1}I^2L^2)^{IL}$ this
commutes with $\alpha=CL^{-2}$. Let us write
$V=C^{-1}I^2=I^{2}C^{-1}$, which will represent the matrix $\begin{pmatrix}0 &-1\\ 1&-1\end{pmatrix}$.\\
 $\alpha^{(IL)^{-1}}=ILCL^{-2}L^{-1}I^{-1}=I^2LI^{-1}=I^2C^{-1}L^{-1}=VL^{-1}$. So
the second commutator after conjugating it by $(IL)^{-1}$ becomes
$[\beta^{-1}\alpha\beta,\alpha]^{(IL)^{-1}}=[(\beta^{-1}\alpha\beta)^{(IL)^{-1}},\alpha^{(IL)^{-1}}]=[L^2C^{-1}I^2L^2,VL^{-1}]
=L^{-2}V^{-1}L^{-2}LV^{-1}L^{2}VL^{2}VL^{-1}$. Notice that
$(CIL)^2=CI(LCI)L=CI(ICL)L=VL^2$, so commutating now by $L^{-2}$ we
get
$V^{-1}L^{-1}V^{-1}L^2VL^2VL^2=V^{-1}L^{-1}V^{-2}(VL^2)^3=CI^2(CL)^{-1}(CIL)^6=(CIL)^7$.\\
 So the proof of the theorem is accomplished.

\bibliographystyle{amsalpha}

\end{document}